\documentclass[a4paper,10pt]{article}
\usepackage{stmaryrd}
\usepackage{amsfonts}
\usepackage{bbm}
\usepackage{amscd}
\usepackage{mathrsfs}
\usepackage{latexsym,amssymb,amsmath,amscd,amscd,amsthm,amsxtra,xypic}
\usepackage[dvips]{graphicx}
\usepackage[utf8]{inputenc}
\usepackage[T1]{fontenc}
\usepackage{lmodern}
\usepackage{amssymb}
\usepackage[all]{xy}
\usepackage{nicefrac,mathtools,enumitem}
\usepackage{microtype}

\newtheorem{thm}{Theorem}[section]
\newtheorem{lem}[thm]{Lemma}
\newtheorem{cor}[thm]{Corollary}
\newtheorem{pro}[thm]{Proposition}
\newtheorem{ex}[thm]{Example}
\newtheorem{rmk}[thm]{Remark}
\newtheorem{defi}[thm]{Definition}

\newcommand{\lon }{\,\rightarrow\,}
\newcommand{\be }{\begin{equation}}
\newcommand{\ee }{\end{equation}}

\newcommand{\pf}{\noindent{\bf Proof.}\ }

\newcommand{\Real}{\mathbb R}

\newcommand{\huaL}{\mathcal{L}}

\newcommand{\huaF}{\mathcal{F}}

\newcommand{\huaV}{\mathcal{V}}

\newcommand{\huaP}{\mathcal{P}}

\newcommand{\huaO}{\mathcal{O}}

\newcommand{\CWM}{C^{\infty}(M)}

\newcommand{\frka}{\mathfrak a}

\newcommand{\frkd}{\mathfrak d}

\newcommand{\frkg}{\mathfrak g}

\newcommand{\frkD}{\mathfrak D}

\newcommand{\frkX}{\mathfrak X}

\def\qed{\hfill ~\vrule height6pt width6pt depth0pt}

\newcommand{\br}[1]{   [ \cdot,    \cdot  ]   }

\newcommand{\dev}{\mathfrak{D}}

\newcommand{\id}{\rm{id}}
\newcommand{\Id}{\mathrm{Id}}

\newcommand{\g}{\mathfrak g}

\newcommand{\dM}{\mathrm{d}}

\newcommand{\Hom}{\mathrm{Hom}}
\newcommand{\Der}{\mathrm{Der}}

\newcommand{\gl}{\mathfrak {gl}}

\newcommand{\End}{\mathrm{End}}
\newcommand{\ad}{\mathrm{ad}}

\newcommand{\Img}{\mathrm{im}}

\begin{document}
\title{
{Left-symmetric algebroids
\thanks
 {
Research partially supported by NSFC  (11001133, 11101179, 11271202, 11221091) and SRFDP
(20100061120096, 20120031110022).
 }
} }
\author{Jiefeng Liu$^1$, Yunhe Sheng$^1$, Chengming Bai$^2$ and Zhiqi Chen$^3$\\
$^1$Department of Mathematics, Jilin University,
 \\Changchun 130012, Jilin, China\\
$^2$Chern Institute of Mathematics and LPMC, Nankai University,\\ Tianjin 300071, China\\
$^3$School of Mathematical Sciences and LMPC, Nankai University, \\Tianjin 300071, China  \\
\quad Email:~liujf12@126.com,\quad shengyh@jlu.edu.cn,\\ baicm@nankai.edu.cn,\quad chenzhiqi@nankai.edu.cn
}

\date{}


\maketitle
\begin{abstract}
In this paper, we introduce a notion of a left-symmetric algebroid,
 which is a generalization of a left-symmetric algebra from a vector space to a vector bundle. The left multiplication gives rise to a representation of the corresponding sub-adjacent Lie algebroid. We construct left-symmetric algebroids from $\mathcal O$-operators on Lie algebroids. We study phase spaces of Lie algebroids in terms of left-symmetric algebroids. Representations of left-symmetric algebroids are studied in detail. At last, we study deformations of left-symmetric algebroids, which could be controlled by the second cohomology class in the deformation cohomology.
\end{abstract}

\tableofcontents
\section{Introduction}

Left-symmetric
algebras (also called  pre-Lie algebras, quasi-associative algebras, Vinberg
algebras and so on) are a class of nonassociative algebras coming
from the study of convex homogeneous cones, affine manifolds and affine structures on
Lie groups, deformation of associative algebras and then appear in many fields in
mathematics and mathematical physics, such as complex and symplectic structures on Lie
groups and Lie algebras, integrable systems, Poisson brackets and infinite dimensional
Lie algebras, vertex algebras, quantum field theory, operads
 and so on.
See \cite{Andrada,Lichnerowicz,Bakalov,pre-Lie operad}, and the survey \cite{Pre-lie algebra in geometry} and the references therein for more details.

  The beauty of a left-symmetric algebra is that the
commutator gives rise to a Lie algebra and the left multiplication
gives rise to a representation of the commutator Lie algebra. So
left-symmetric algebras naturally play important roles in the study
involving the representations of Lie algebras  on the underlying
spaces of the Lie algebras themselves or their dual spaces. For
example, they are the underlying algebraic structures of the
non-abelian phase spaces of Lie algebras
\cite{Kupershmidt1,non-abelian phase spaces}, which lead to a
bialgebra theory of left-symmetric algebras \cite{Left-symmetric
bialgebras}. They are also regarded as the algebraic structures
``behind'' the classical Yang-Baxter equations and they provide a
construction  of solutions of classical Yang-Baxter equations in
certain semidirect product Lie algebra structures (that is, over the
``double'' spaces) induced by left-symmetric algebras
\cite{Kupershmidt2,Bai:CYBE}.

The notion of a Lie algebroid was introduced by Pradines in 1967, which
is a generalization of Lie algebras and tangent bundles. A Lie algebroid structure on a vector bundle $A\longrightarrow M$ is
a pair that consists of a Lie algebra structure $[\cdot,\cdot]_A$ on
the section space $\Gamma(A)$ and a vector bundle morphism
$a_A:A\longrightarrow TM$, called the anchor, such that the
following relations are satisfied:
$$~[x,fy]_A=f[x,y]_A+a_A(x)(f)y,\quad \forall~x,y\in\Gamma(A),~f\in
\CWM.$$
We usually denote a Lie algebroid by $(A,[\cdot,\cdot]_A,a_A)$, or $A$
if there is no confusion. See \cite{General theory of Lie groupoid and Lie algebroid} for more
details about Lie algebroids.

In this paper, we introduce a notion of a left-symmetric algebroid, which is a geometric generalization of a left-symmetric algebra, such that we have the following diagram:
$$
\xymatrix{
 \mbox{Lie~ algebras} \ar[d]_{\mbox{geometric~ generalization}}
                && \mbox{left-symmetric  ~algebras} \ar[d]^{\mbox{geometric~ generalization}}\ar[ll]_{\mbox{commutator}\qquad}  \\
\mbox{Lie~ algebroids}
                && \mbox{left-symmetric  ~algebroids.}    \ar[ll]_{\mbox{commutator}\qquad}       }
$$
We study phase spaces of Lie algebroids in terms of left-symmetric
algebroids. In particular, we construct a paracomplex structure and
a complex structure on the phase space of the sub-adjacent Lie
algebroid associated to a left-symmetric algebroid. We study
representations of left-symmetric algebroids. Different from the
case of a left-symmetric algebra, there is no regular representation
associated to a left-symmetric algebroid. This forces us to
introduce a new cohomology, which is called the deformation
cohomology, when we study  deformations of a left-symmetric
algebroid.

The paper is organized as follows. In Section 2, we recall notions of left-symmetric algebras, Lie-admissible algebras and representations of Lie algebroids. In Section 3, we study basic properties of a left-symmetric algebroid. Similar as the case of a left-symmetric algebra, we can obtain the sub-adjacent Lie algebroid from a left-symmetric algebroid by using the commutator. The left multiplication gives rise to a representation of the sub-adjacent Lie algebroid (See Theorem \ref{thm:sub-adjacent}). We construct left-symmetric algebroids using $\mathcal O$-operators. We extend the multiplication on $\Gamma(A)$ to $\Gamma(\Lambda^\bullet A)$, and obtain a graded Lie-admissible algebra (see Theorem \ref{thm:graded Lie-adm}). In Section 4, we study phase spaces of Lie algebroids in terms of left-symmetric algebroids. There is a phase space $\huaP$ of the sub-adjacent Lie algebroid associated to a left-symmetric algebroid naturally (see Theorem \ref{thm:phase space}). Furthermore, there is a natural paracomplex structure on the phase space $\huaP$ (see Proposition \ref{pro:paracomplex}). The notion of a quadratic left-symmetric algebroid is introduced, and there is a natural complex structure on the phase space $\huaP$ of the   sub-adjacent Lie algebroid associated to a quadratic left-symmetric algebroid (see Proposition \ref{pro:complex}). In Section 5, we study representations of left-symmetric algebroids. In the last Section, we introduce the deformation cohomology associated to a left-symmetric algebroid, which could control deformations. We also introduce a notion of a Nijenhuis operator, which could generate a trivial deformation (see Theorem \ref{thm:trivial def}).

In this paper, all the vector bundles are over the same manifold
$M$. For two vector bundles $A$ and $B$, a bundle map from $A$ to
$B$ is a base-preserving map and $\CWM$-linear.

\section{Preliminaries}

In this section, we briefly recall left-symmetric algebras, Lie-admissible algebras and  representations of Lie algebroids.

\begin{defi}  A {\bf left-symmetric algebra} is a pair $(\frkg,\cdot_\frkg)$, where $\g$ is a vector space, and  $\cdot_\frkg:\g\otimes \g\longrightarrow\g$ is a bilinear multiplication
satisfying that for any $x,y,z\in \g$, the associator
$(x,y,z)=(x\cdot_\frkg y)\cdot_\frkg z-x\cdot_\frkg(y\cdot_\frkg z)$ is symmetric in $x,y$,
i.e.,
$$(x,y,z)=(y,x,z),\;\;{\rm or}\;\;{\rm
equivalently,}\;\;(x\cdot_\frkg y)\cdot_\frkg z-x\cdot_\frkg(y\cdot_\frkg z)=(y\cdot_\frkg x)\cdot_\frkg
z-y\cdot_\frkg(x\cdot_\frkg z).$$
\end{defi}
 For any $x\in \g$, let $L_x$
denote the left multiplication operator, i.e. $L_x(y)=x\cdot_\frkg
y$ for any $y\in \g$. The following conclusion is known
\cite{Pre-lie algebra in geometry}:
\begin{lem} \label{lem:sub-ad} Let $(\g,\cdot_\g)$ be a left-symmetric algebra. The commutator
$ [x,y]_\g=x\cdot_\g y-y\cdot_\g x$ defines a Lie algebra $G(\g)$,
which is called the {\bf sub-adjacent Lie algebra} of $\g$ and $\g$ is also
called a {\bf compatible left-symmetric algebra}  on the Lie
algebra $G(\g)$. Furthermore, $L:G(\g)\rightarrow
\gl(\g)$ with $x\rightarrow L_x$ gives a representation of the Lie
algebra $G(\g)$.
\end{lem}

A {\bf Lie-admissible algebra} is a nonassociative algebra $(\g,\cdot_\g)$ whose commutator algebra is a Lie algebra.
More precisely, it is equivalent to the following equation:
\begin{equation}
  (x,y,z)-(y,x,z)+(y,z,x)-(z,y,x)+(z,x,y)-(x,z,y)=0,\;\;\forall x,y,z\in \frak g.
\end{equation}

Obviously, a left-symmetric algebra is a Lie-admissible algebra.

For a vector bundle $E\longrightarrow M$, its gauge Lie algebroid
$\frkD(E)$ is just the gauge Lie algebroid of the
 frame bundle
 $\huaF(E)$, which is also called the covariant differential operator bundle of $E$.   Here we treat each
element $\frkd$ of $\frkD(E)$ at $m\in M$ as an $\Real$-linear
operator $\Gamma(E)\lon E_m$ together with some $x\in T_mM$ (which
is uniquely determined by $\frkd$ and called the anchor of
$\frkd$) such that
$$
\frkd(fu)=f(m)\frkd(u)+x(f)u(m), \, \quad \quad ~~ \forall~
f\in\CWM,
 u\in\Gamma(E).
$$
It is  known that $\frkD(E)$ is a  transitive Lie algebroid over $M$.
The anchor of $\dev{E}$ is given by $a(\frkd)=x$ and the Lie bracket
$[\cdot,\cdot ]_\frkD$ of $\Gamma(\frkD(E))$ is given by the usual
commutator of two operators.

Let $(A,[\cdot,\cdot]_A,a_A),(B,[\cdot,\cdot]_B,a_B)$ be two Lie
algebroids (with the same base), a {\bf base-preserving morphism}
from $A$ to $B$ is a bundle map $\sigma:A\longrightarrow B$ such
that
\begin{eqnarray*}
  a_B\circ\sigma=a_A,\quad
  \sigma[x,y]_A=[\sigma(x),\sigma(y)]_B.
\end{eqnarray*}

A {\bf representation} of a  Lie algebroid $(A,[\cdot,\cdot]_A,a_A)$
 on $E$ is a base-preserving morphism $\rho$ from $A$ to the Lie algebroid $\frkD(E)$.
Denote a representation by $(E;\rho).$ The dual representation of a
Lie algebroid $A$ on $E^*$  is the bundle map
$\rho^*:A\longrightarrow \frkD(E^*)$ given by
$$
\langle \rho^*(x)(\xi),y\rangle=a_A(x)\langle \xi,y\rangle-\langle
\xi,\rho(x)(y)\rangle,\quad \forall~x,y\in
\Gamma(A),~\xi\in\Gamma(E^*).
$$

A Lie algebroid $(A,[\cdot,\cdot]_A,a_A)$ naturally represents on the trivial line bundle $E=M\times \mathbb R$ via the anchor map $a_A:A\longrightarrow TM$. The corresponding coboundary operator $d:\Gamma(\Lambda^kA^*)\longrightarrow \Gamma(\Lambda^{k+1}A^*)$ is given by
\begin{eqnarray*}
  d\varpi(x_1,\cdots,x_{k+1})&=&\sum_{i=1}^{k+1}(-1)^{i+1}a_A(x_i)\varpi(x_1\cdots,\widehat{x_i},\cdots,x_{k+1})\\
  &&+\sum_{i<j}(-1)^{i+j}\varpi([x_i,x_j]_A,x_1\cdots,\widehat{x_i},\cdots,\widehat{x_j},\cdots,x_{k+1}).
\end{eqnarray*}
In particular, a $2$-form $\varpi\in\Gamma(\Lambda^2A^*)$ is a {\bf 2-cocycle} if $d\varpi=0$, i.e.
\begin{equation}
  a_A(x)\varpi(y,z)- a_A(y)\varpi(z,x)+ a_A(z)\varpi(x,y)-\varpi([x,y]_A,z)+\varpi([x,z]_A,y)-\varpi([y,z]_A,x).
\end{equation}

\section{Some basic properties of a left-symmetric algebroid}

\begin{defi}\label{defi:left-symmetric algebroid}
A {\bf left-symmetric algebroid} structure on a vector bundle
$A\longrightarrow M$ is a pair that consists of a left-symmetric
algebra structure $\cdot_A$ on the section space $\Gamma(A)$ and a
vector bundle morphism $a_A:A\longrightarrow TM$, called the anchor,
such that for all $f\in\CWM$ and $X,Y\in\Gamma(A)$, the following
relations are satisfied:
\begin{itemize}
\item[\rm(i)]$~x\cdot_A(fy)=f(x\cdot_A y)+a_A(x)(f)y,$
\item[\rm(ii)] $(fx)\cdot_A y=f(x\cdot_A y).$
\end{itemize}
\end{defi}

We usually denote a left-symmetric algebroid by $(A,\cdot_A, a_A)$. A left-symmetric algebroid $(A,\cdot_A,a_A)$ is called {\bf transitive (regular)} if  $a_A$ is
 surjective (of constant rank).
Any left-symmetric  algebra is a left-symmetric algebroid. Let $M$
be a differential manifold with a flat torsion free connection
$\nabla$, then $(TM,\nabla,\id)$ is a left-symmetric algebroid.

For any $x\in\Gamma(A)$, we can define
$L_x:\Gamma(A)\longrightarrow\Gamma(A)$ by $L_x(y)=x\cdot_A y$.
Condition (i) in the above definition means that $L_x\in \frkD(A)$.
Condition (ii) means that the map $x\longmapsto L_x$ is
$C^\infty$-linear. Thus, $L:A\longrightarrow \frkD(A)$ is a bundle
map.

\begin{thm}\label{thm:sub-adjacent}
  Let $(A,\cdot_A, a_A)$ be a left-symmetric algebroid. Define  a skew-symmetric bilinear bracket operation $[\cdot,\cdot]_A$ on $\Gamma(A)$ by
  $$
  [x,y]_A=x\cdot_A y-y\cdot_A x,\quad \forall ~x,y\in\Gamma(A).
  $$
Then, $(A,[\cdot,\cdot]_A,a_A)$ is a Lie algebroid, and denoted by
$G(A)$, called the {\bf sub-adjacent Lie algebroid} of
 $(A,\cdot_A,a_A)$. Furthermore, $L:A\longrightarrow \frkD(A)$ with
$x\longmapsto L_x$ gives a
  representation of the Lie algebroid  $(A,[\cdot,\cdot]_A,a_A)$.
\end{thm}
\pf Since $(\Gamma(A),\cdot_A)$ is a left-symmetric  algebra,  we
have that $(\Gamma(A),[\cdot,\cdot]_A)$ is a Lie algebra. For any
$f\in \CWM$, by direct computations, we have
\begin{eqnarray*}
  [x,fy]_A&=&x\cdot_A(fy)-(fy)\cdot_A x=f(x\cdot_A y)+a_A(x)(f)y-f(y\cdot_A x)\\
  &=&f[x,y]_A+a_A(x)(f)y,
\end{eqnarray*}
which implies that $(A,[\cdot,\cdot]_A,a_A)$ is a Lie algebroid.

To see that the bundle map $L:A\longrightarrow \frkD(A)$ is a
representation, we only need to show that
$L_{[x,y]_A}=[L_x,L_y]_\frkD$, which follows directly from the fact
that $(\Gamma(A),\cdot_A)$ is a left-symmetric  algebra. This
finishes the proof.\qed

\begin{defi}
Let $(\frkg,\cdot_{\frkg})$ be a left-symmetric  algebra. An {\bf action}
of $\frkg$ on $M$ is a linear map $\rho:\frkg\longrightarrow\frkX(M)$  from $\g$ to the space of vector fields on $M$,
such that for all $x,y\in\frkg$, we have
 $$
 \rho(x\cdot_{\frkg}y-y\cdot_{\frkg}x)=[\rho(x),\rho(y)]_{\frkX(M)}.
 $$
\end{defi}
It is obvious that $\rho$ is also an action of the corresponding
sub-adjacent Lie algebra $G(\frkg)$ on $M$.

 Given an action of
$\frkg$ on $M$, let $A=M\times
\g$ be the trivial bundle. Define an anchor map $a_\rho:A\longrightarrow TM$ and a multiplication $\cdot_\rho:\Gamma(A)\times \Gamma(A)\longrightarrow \Gamma(A)$   by
\begin{eqnarray}
a_\rho(m,u)&=&\rho(u)_m,\quad \forall ~m\in M, u\in\g,\label{action pre2}\\
x\cdot_\rho y&=&\huaL_{\rho(x)}(y)+x\cdot_\g y, \quad \forall~x,y\in\Gamma(A),\label{action pre1}
\end{eqnarray}
where  $x\cdot_\g y$ is the pointwise multiplication.
\begin{pro}
With the above notations, $(A=M\times\frkg,\cdot_\rho,a_\rho)$ is a
left-symmetric algebroid, which we call  an {\bf action left-symmetric algebroid},
where $\cdot_\rho $ and $a_\rho$ are given by $(\ref{action pre1})$
and $(\ref{action pre2})$ respectively.
\end{pro}
\pf For all $u,v,w \in
\frkg$ and $f,g,h \in C^\infty(M)$, we have
\begin{eqnarray*}
&&((fu)\cdot_\rho(gv))\cdot_\rho(hw)-(fu)\cdot_\rho((gv)\cdot_\rho(hw))
-((gv)\cdot_\rho(fu))\cdot_\rho(hw)-(gv)\cdot_\rho((fu)\cdot_\rho(hw))\\
&=&-fg\rho(u)\rho(v)(h)w+fg\rho(v)\rho(u)(h)w+fg\rho(u\cdot_A
v)(h)w-fg\rho(v\cdot_A u)(h)w\\
&=& -fg[\rho(u),\rho(v)]_{\frkX(M)}(h)w+fg\rho([u,v]_A)(h)w=0,
\end{eqnarray*}
which implies that $(\Gamma(A),\cdot_\rho)$ is a left-symmetric  algebra.
Furthermore, we have
\begin{eqnarray*}
x\cdot_\rho (fy)&=&\huaL_{\rho(x)}(fy)+x\cdot_\g (fy)\\
&=&f\huaL_{\rho(x)}(y)+\rho(x)(f)y+f(x\cdot_\g y)\\
&=&f(x\cdot_\rho(fy))+a_\rho(x)(f)y;\\
(fx)\cdot_\rho y&=&\huaL_{\rho(fx)}(y)+(fx)\cdot_\g y\\
&=&f(x\cdot_\rho y).
\end{eqnarray*}
Thus, $(A,\cdot_\rho,a_\rho)$ is a left-symmetric algebroid. \qed\vspace{3mm}

Obviously, the sub-adjacent Lie
algebroid of an action left-symmetric algebroid is an action Lie algebroid. See \cite[Example 3.3.7]{General theory of Lie groupoid and Lie algebroid} for more details about action Lie algebroids.

\begin{defi}
Let $(A_1,\cdot_1,a_1)$ and $(A_2,\cdot_2,a_2)$ be left-symmetric
algebroids on $M$, a bundle map $\varphi:A_1\longrightarrow A_2$ is
called a {\bf homomorphism}  of left-symmetric algebroids, if the following
conditions are satisfied:
$$\varphi(x \cdot_1 y)=\varphi(x)\cdot_2\varphi(y),\quad
a_2\circ\varphi=a_1, \quad\forall x,y\in \Gamma(A_1). $$
\end{defi}

It is straightforward to obtain following proposition.
\begin{pro}
Let $(A_1,\cdot_1,a_1)$ and $(A_2,\cdot_2,a_2)$ be left-symmetric
algebroids, and  $\varphi:A_1\longrightarrow A_2$ a homomorphism of
left-symmetric algebroids. Then, $\varphi$ is also a Lie algebroid
homomorphism of the corresponding sub-adjacent Lie algebroids.
\end{pro}

\begin{defi}
  Let $(A,[\cdot,\cdot]_A,a_A)$ be a Lie algebroid and $\rho:A\longrightarrow \frkD(E)$ be its  representation, where $E$ is a vector bundle over the same base manifold $M$. A bundle map
  $T:E\longrightarrow A$ is called an {\bf $\mathcal{O}$-operator} if
  $$
  [T(u),T(v)]_A=T(\rho(T(u))(v)-\rho(T(v))(u)), \quad\forall~u,v\in\Gamma(E).
  $$
\end{defi}

\begin{rmk}
  Consider the semidirect product Lie algebroid $(A\ltimes_{\rho} E,[\cdot,\cdot]_s,\frka)$, where $\frka(x+u)=a_A(x)$ and the bracket $[\cdot,\cdot]_s$ is given by
  $$
  [x+u,y+v]_s=[x,y]_A+\rho(x)(v)-\rho(y)(u).
  $$ For any $\huaO$-operator $T:E\longrightarrow A$, $\tilde{T}=\left(\begin{array}{cc}
    0&T\\0&0
  \end{array}\right)$ is a Nijenhuis operator for the Lie algebroid $A\ltimes_\rho
  E$. More precisely, we have
  $$
  [\tilde{T}(x+u),\tilde{T}(y+v)]_s=\tilde{T}([\tilde{T}(x+u),y+v]_s+[x+u,\tilde{T}(y+v)]_s-\tilde{T}[x+u,y+v]_s).
  $$
  See \cite{Dorfman} for more details about Nijenhuis operators and their applications.
\end{rmk}

Let $T:E\longrightarrow A$ be an $\huaO$-operator, define a bilinear
operator $\cdot_E$ on $\Gamma(E)$ by
$$
u\cdot_E v=\rho(T(u))(v).
$$

\begin{pro}
  With the above notations, $(E,\cdot_E,a_E=a_A\circ T)$ is a left-symmetric algebroid, and $T$ is a base-preserving Lie algebroid homomorphism from the sub-adjacent Lie algebroid to $(A,[\cdot,\cdot]_A,a_A)$.
\end{pro}
\pf It is not hard to see that $(\Gamma(E),\cdot_E)$ is a left-symmetric
algebra. For any $f\in\CWM$, we have
\begin{eqnarray*}
  (fu)\cdot_E v&=&\rho(T(fu))(v)=\rho(fT(u))(v)=f\rho(T(u))(v),\\
  u\cdot_E(fv)&=&\rho(T(u))(fv)=f\rho(T(u))(v)+a_A\circ T(u)(f)v.
\end{eqnarray*}
Thus, $(E,\cdot_E,a_A\circ T)$ is a left-symmetric algebroid.

Let $[\cdot,\cdot]_E$ be the sub-adjacent Lie bracket on $\Gamma(E)$. Then
we have
$$
T[u,v]_E=T(u\cdot_E v-v\cdot_E
u)=T(\rho(T(u))(v)-\rho(T(v))(u))=[T(u),T(v)]_A.
$$
So $T$ is a base-preserving homomorphism between Lie
algebroids.\qed\vspace{3mm}

Given a left-symmetric algebroid $(A,\cdot_A,a_A)$, denote by $\huaV^k(A)=\Gamma(\Lambda^kA)$. We extend the multiplication $\cdot_A$ to
$\huaV(A):=\oplus_{k}\huaV^{k}(A)$, for which we use the notation $\cdot_S$, by the following rules:

\begin{itemize}
\item[$\rm(a)$]$\huaV^{k}(A)\cdot_S\huaV^{l}(A)\subset\huaV^{k+l-1}(A)$ for all
$k$ and $l$;
\item[$\rm(b)$]for any smooth function $f\in\huaV^{0}=C^\infty(M)$
and $x\in\huaV^{1}(A)$,
\begin{eqnarray*}
x\cdot_S f=a_A(x)(f),\quad f\cdot_S x=0;
\end{eqnarray*}
\item[$\rm(c)$]for any $x,y\in\huaV^{1}(A)$, $x\cdot_S y=x\cdot_A y;$
\item[$\rm(d)$]for $x$ and $y$ of degrees $|x|$ and $|y|$
respectively,
\begin{eqnarray*}
x\cdot_S(y\wedge z)=(x\cdot_S y)\wedge z+(-1)^{(|x|-1)|y|}y\wedge
(x\cdot_S
z),\\
(x\wedge y)\cdot_S z=x\wedge(y\cdot_S z)+(-1)^{(|z|-1)|y|}(x\cdot_S
z)\wedge y.
\end{eqnarray*}
\end{itemize}

Through a careful calculation, we have
\begin{lem} For any $x_1,\cdots,x_k,y_1,\cdots, y_l\in\Gamma(A)$, we have
\begin{eqnarray}
&&(x_1\wedge x_2\cdots\wedge x_k)\cdot_S(y_1\wedge y_2\cdots\wedge
y_l)\nonumber\\
&=&\sum_{i,j}(-1)^{i+j}x_i\cdot_A y_j\wedge x_1\cdots\wedge
\hat{x_i}\wedge\cdots \wedge x_k\wedge y_1\cdots
\hat{y_j}\wedge\cdots y_l.\label{S operation}
\end{eqnarray}
\end{lem}

Now we introduction some notations which will be used below.
\begin{eqnarray}
\sigma(x): &=&|x|-1;\nonumber \\
{[x,y]_S} &=&x\cdot_S y-(-1)^{(|x|-1)(|y|-1)}y\cdot_S x;\label{S bracket}\\
C(x,y,z): &=& (x\cdot_S y)\cdot_S z-x\cdot_S(y\cdot_S z);\nonumber\\
{[x,y,z]_S}: &=& C(x,y,z)-(-1)^{(|x|-1)(|y|-1)}C(y,x,z);\label{left-symmetry super}\\
CI(x,y,z):&=&(-1)^{(|x|-1)(|z|-1)}[x,y,z]_S+(-1)^{(|y|-1)(|x|-1)}[y,z,x]_S\nonumber\\
&&+(-1)^{(|z|-1)(|y|-1)}[z,x,y]_S.\label{jie hezi}
\end{eqnarray}
It is obvious that
\begin{eqnarray*}
 \sigma([x,y,z]_S)&=&\sigma(x)+\sigma(y)+\sigma(z),\\
~[x,y,z]_S&=&-(-1)^{\sigma(x)\sigma(y)}[y,x,z]_S.
\end{eqnarray*}

 Define
$\tilde{\huaV}^{k}(A):=\huaV^{k+1}(A)$. In general, $(\tilde{\huaV}(A):=
\oplus_{k}\tilde{\huaV}^{k}(A),\cdot_S)$ is not a graded left-symmetric algebra. But, we have

\begin{thm}\label{thm:graded Lie-adm}
With the above notations, $(\tilde{\huaV}(A),\cdot_S)$ is a graded
Lie-admissible algebra. Furthermore, the graded Lie bracket
$[\cdot,\cdot]_S$ given by \eqref{S bracket} is exactly the
Schouten bracket on  $\huaV(A)$ corresponding to the sub-adjacent
Lie algebroid $(A,[\cdot,\cdot]_A,a_A)$.
\end{thm}
\pf By $(\ref{S operation})$, we can obtain $CI=0$ directly, which
implies that $(\tilde{\huaV}(A),\cdot_S)$ is a graded Lie-admissible
algebra. Since for any $x,y\in\Gamma(A)$, $[x,y]_S=[x,y]_A$, to see
that $[\cdot,\cdot]_S$ is the Schouten bracket, we only need to
prove that $[\cdot,\cdot]_S$ satisfies the graded Leibniz rule. For
any $x\in\huaV^{|x|}(A)$, $y\in\huaV^{|y|}(A)$ and
$z\in\huaV^{|z|}(A)$, we have
\begin{eqnarray*}
[x,y\wedge z]_S&=&x\cdot_S (y\wedge
z)-(-1)^{(|x|-1)(|y|+|z|-1)}(y\wedge
z)\cdot_S x\\
&=& (x\cdot_S y)\wedge z+(-1)^{(|x|-1)|y|}y\wedge (x\cdot_S
z)-(-1)^{(|x|-1)(|y|+|z|-1)}\\
&&(y\wedge(z\cdot_S x)+(-1)^{(|x|-1)|z|}(y\cdot_S x)\wedge z) \\&=&
[x,y]_S\wedge z+(-1)^{(|x|-1)|y|}y\wedge[x,z]_S,
\end{eqnarray*}
which finishes the proof.\qed

\section{Phase spaces of Lie algebroids}

On the direct sum vector bundle $A\oplus A^*$, there is a natural $2$-form $\omega\in \Lambda^2\Gamma(A\oplus A^*)$ given by
\begin{equation}\label{eq:o}
\omega(x+\xi, y+\eta)=\langle\eta,x\rangle-\langle\xi,y\rangle,\quad
\forall~x,y\in\Gamma(A), \xi,\eta\in\Gamma(A^*).
\end{equation}


\begin{defi}
Let $(A,[\cdot,\cdot]_A,a_A)$ be a Lie algebroid. Suppose that there
is a Lie algebroid structure on the dual bundle  $A^*$. If there is
a Lie algebroid structure on $\huaP=A\oplus{A^*}$ such that $A$ and
$A^*$ are subalgebroids and the $2$-form $\omega$ given by
\eqref{eq:o} is a $2$-cocycle, then $\huaP$ is called a {\bf phase
space} of the Lie algebroid $A$.
\end{defi}

Given a Lie algebroid  $(A,[\cdot,\cdot]_A,a_A)$ and a representation $\rho:A\longrightarrow \frkD(A)$. Let $\rho^*:A\longrightarrow \frkD(A^*)$ be its
  dual representation. Then we have the semidirect product Lie algebroid $A\ltimes_{\rho^*}A^*$.  In particular, let $(A,\cdot_A,a_A)$ be a left-symmetric algebroid. Then we consider the
semidirect product Lie algebroid $G(A)\ltimes_{L^*} A^*$.

\begin{thm}\label{thm:phase space}
  With the above notations, $\omega$ which is given by \eqref{eq:o}, is a nondegenerate closed $2$-form on the Lie algebroid $G(A)\ltimes_{L^*} A^*$. Thus, $G(A)\ltimes_{L^*} A^*$ is a phase space of the Lie algebroid $G(A)$.
    \end{thm}
 \pf   It is obvious that $\omega$ is nondegenerate.  Next we show that $\omega$ is closed, i.e. $d\omega=0$. For all $x,y,z\in\Gamma(A)$ and $\xi,\eta,\gamma\in\Gamma(A^*)$, we have
  \begin{eqnarray*}
    &&d\omega(x+\xi, y+\eta,z+\gamma)\\&=&a_A(x)\omega(y+\eta,z+\gamma)+c.p.
    -\omega([x+\xi, y+\eta]_s,z+\gamma)+c.p.\\
    &=&a_A(x)(\langle\gamma,y\rangle-\langle\eta,z\rangle)+c.p.-\omega(x\cdot_A y-y\cdot_A x+L_x^*\eta-L_y^*\xi,z+\gamma)+c.p.\\
    &=&a_A(X)(\langle\gamma,Y\rangle-\langle\eta,Z\rangle)+c.p.\\
    &&-\langle x\cdot_A y,\gamma\rangle+\langle y\cdot_A x,\gamma\rangle-\langle\eta, x\cdot_A z\rangle+a_A(x)\langle\eta,z\rangle-\langle\xi, y\cdot_A z\rangle+a_A(y)\langle\xi,z\rangle+c.p.\\
    &=&0,
  \end{eqnarray*}
  which finishes the proof.\qed\vspace{3mm}



Conversely, we have

\begin{thm}\label{thm:phasespace}
Let $(A,[\cdot,\cdot]_A,a_A)$ be a Lie algebroid and $\rho:A\longrightarrow \frkD(A)$ be a representation.
If~ $\huaP=A\ltimes_{\rho^*}A^*$ is a phase of a Lie algebroid $(A,[\cdot,\cdot]_A,a_A)$, then $(A,\cdot_A,a_A)$ is a left-symmetric algebroid  whose sub-adjacent Lie algebroid is exactly
  $(A,[\cdot,\cdot]_A,a_A)$, where $\cdot_A$ is given by $$x\cdot_A y=\rho(x)(y),\quad\forall x,y\in \Gamma(A). $$
 Furthermore, there is a compatible left-symmetric algebroid structure
$(\ast_\huaP,a_\huaP)$ on $\huaP$ given as follows:
\begin{eqnarray*}x\ast_\huaP y&=&x\cdot_A y,\quad \xi\ast_\huaP\eta=0,\quad
x\ast_\huaP\xi=L^*_x\xi,\quad \xi\ast_\huaP x=0,\\
a_\huaP(x+\xi)&=&a_A(x).\end{eqnarray*}
\end{thm}

\pf Since $d\omega=0$,  we have
\begin{eqnarray*}
d\omega(x,y,\gamma)&=&a_A(x)\omega(y,\gamma)-a_A(y)\omega(x,\gamma)-\omega([x,y]_A,\gamma)+\omega([x,\gamma]_s,y)-\omega([y,\gamma]_s,x)\\
&=&a_A(x)\langle y,\gamma\rangle-a_A(y)\langle
x,\gamma\rangle-\langle\gamma,[x,y]_A\rangle-\langle
\rho^*(x)\gamma,y\rangle+\langle
\rho^*(y)\gamma,x\rangle\\
 &=& -\langle \gamma,[x,y]_A-\rho(x)y+\rho(y)x\rangle=0,
\end{eqnarray*}
for all $x,y\in A,~ \xi,\eta,\gamma\in A^*$, that is
\begin{eqnarray}\label{relation Lie and pre-Lie}
 [x,y]_A=\rho(x)y-\rho(y)x.
\end{eqnarray}
It is obvious that $(\Gamma(A),\cdot_A)$ is a left-symmetric
algebra. Furthermore, we have
\begin{eqnarray*}
x\cdot_A (f{y})&=&\rho(x)(f{y})=f\rho(x)y+a_A(x)(f)y,\\
(f{x})\cdot_A y&=&\rho(f{x})y=f\rho(x)y=f (x\cdot_A y),
\end{eqnarray*}
for all $x,y\in \Gamma(A), f\in C^{\infty}(M)$. Consequently,
$(A,\cdot_A,a_A)$ is a left-symmetric algebroid. It is
straightforward to obtain the other conclusions. The proof is
completed. \qed

\begin{defi}
Let $\huaP_1$ and $\huaP_2$ be phase spaces of Lie algebroids $A_1$ and
$A_2$ respectively. $\huaP_1$ is said to be {\bf isomorphic} to $\huaP_2$ if there
exists a Lie algebroid isomorphism $\varphi:\huaP_1\longrightarrow \huaP_2$
satisfying the following conditions:
$$\varphi(A_1)=A_2,\quad \varphi({A^*_1})={A^*_2},\quad \omega_1(X,Y)=\omega_2(\varphi(X),\varphi(Y)),\quad \forall~ X,Y\in \Gamma(\huaP_1),$$
where $\omega_1$ and $\omega_2$ are the natural skew-symmetric $2$-forms on $\huaP_1$ and $\huaP_2$ respectively.
\end{defi}

\begin{pro}
Let $(A_1,\cdot_1,a_1)$ and $(A_2,\cdot_2,a_2)$ be two
left-symmetric algebroids. Then the phase space
$\huaP_1=G(A_1)\ltimes_{L^*}A_1^*$ is isomorphic to the phase
space $\huaP_2=G(A_2)\ltimes_{L^*}A_2^*$ if and only if $A_1$ is
isomorphic to $A_2$ as left-symmetric algebroids.
\end{pro}
\pf Assume that $\varphi:\huaP_1\longrightarrow \huaP_2$ is the
isomorphism of phase spaces. Then, for all $x,y\in
\Gamma(A_1),~\xi\in \Gamma(A_1^*)$, we have
\begin{eqnarray*}
\langle \varphi(x\cdot_1 y),\varphi(\xi)\rangle&=&\omega(
\varphi(x\cdot_1 y),\varphi(\xi))=\omega(x\cdot_1 y,\xi)=\langle x\cdot_1 y,\xi\rangle\\
&=&a_1(x)\langle y,\xi\rangle-\langle y,L^*_x\xi\rangle\\
&=&a_1(x)\omega(y,\xi)-\omega(y,L^*_x\xi)\\
&=&a_1(x)\omega(\varphi(y),\varphi(\xi))-\omega(\varphi(y),\varphi(L_x^*\xi))\\
&=&a_1(x)\langle\varphi(y),\varphi(\xi)\rangle-\langle\varphi(y),L^*_{\varphi(x)}\varphi(\xi)\rangle\\
&=&a_2\circ\varphi(x)\langle\varphi(y),\varphi(\xi)\rangle-\langle\varphi(y),L^*_{\varphi(x)}\varphi(\xi)\rangle\\
&=&\langle \varphi(x)\cdot_2\varphi(y),\varphi(\xi)\rangle.
\end{eqnarray*}
Since $\varphi({A^*_1})={A^*_2}$,  we have $\varphi(x\cdot_1
y)=\varphi(x)\cdot_2\varphi(y)$, which implies that $\varphi|_{A_1}$ is an
isomorphism of left-symmetric algebroids.

Conversely, we will use the compatible left-symmetric algebroid structures
on $(\huaP_1,\ast_1,\frka_1)$ and $(\huaP_2,\ast_2,\frka_2)$ given in Theorem \ref{thm:phasespace} to prove this
fact. Let $\varphi:A_1\longrightarrow A_2$ be an isomorphism of left-symmetric algebroids. Then $\varphi^*:A_2^*\longrightarrow A_1^*$ is invertible and satisfies
$$\langle (\varphi^*)^{-1}(\xi),\varphi(x)\rangle=\langle \xi,x\rangle.$$
Consider the bundle map $\Phi=(\varphi,(\varphi^*)^{-1}):\huaP_1\longrightarrow \huaP_2$.
For any $x,y\in A_1,~\xi\in A_1^*$, we have
\begin{eqnarray*}
\langle(\varphi^*)^{-1}(x\ast_1\xi),\varphi(y)\rangle&=&\langle
x\ast_1\xi,y\rangle=a_1(x)\langle\xi,y\rangle-\langle\xi,x\cdot_1 y\rangle\\
&=&a_1(x)\langle(\varphi^*)^{-1}(\xi),\varphi(y)\rangle-\langle(\varphi^*)^{-1}(\xi),\varphi(x\cdot_1
y)\rangle \\
&=&a_2\circ\varphi(x)\langle(\varphi^*)^{-1}(\xi),\varphi(y)\rangle-\langle(\varphi^*)^{-1}(\xi),\varphi(x)\cdot_2
\varphi(y)\rangle \\
&=&\langle\varphi(x)\ast_2(\varphi^*)^{-1}(\xi),\varphi(y)\rangle.
\end{eqnarray*}
So we have that
$(\varphi^*)^{-1}(x\ast_1\xi)=\varphi(x)\ast_2(\varphi^*)^{-1}(\xi)$,
i.e. $\Phi[x,\xi]=[\Phi(x),\Phi(\xi)]$.

Furthermore, we have
\begin{eqnarray*}
\omega(x+\xi,y+\eta)&=&\langle x,\eta\rangle-\langle y,\xi\rangle\\
&=& \langle \varphi(x),\varphi^{-1}(\eta)\rangle-\langle
\varphi(y),\varphi^{-1}(\xi)\rangle =\omega(\Phi(x+\xi),\Phi(y+\eta)).
\end{eqnarray*}
Thus, $\Phi$ is an isomorphism of phase spaces. \qed

\begin{defi}
Let $(A,[\cdot,\cdot]_A,a_A)$ be a real Lie algebroid on $M$, a {\bf
paracomplex  structure} on $A$ is a bundle map  $P:A\longrightarrow
A$, satisfying $P^2=\id$
 and the integrability  condition:
$$P[x,y]_A=[P(x),y]_A+[x,P(y)]_A-P[P(x),P(y)]_A,\quad\forall x,y\in\Gamma(A).$$
\end{defi}

\begin{pro}\label{pro:paracomplex}
Let $(A,\cdot_A,a_A)$ be a left-symmetric algebroid. Then on the phase space
$\huaP=G(A)\ltimes_{L^*}A^*$, there is a paracomplex structure
$P:\huaP\longrightarrow \huaP$ given by $\left(\begin{array}{cc}
    \id&0\\0&{-\id}
  \end{array}\right)$.
  More precisely, we have
\begin{eqnarray}\label{paracomplex structure}
P(x+\xi)=x-\xi,\quad\forall x\in \Gamma(A), ~\xi\in\Gamma(A^*).
\end{eqnarray}
\end{pro}
\pf It is obvious that $P^2=\id$. For any $x,y\in \Gamma(A)$ and
$\xi,\eta\in \Gamma(A^*)$, we have
\begin{eqnarray*}
P[x+\xi,y+\eta]_s &=&P([x,y]_A+L^*_x\eta-L^*_y\xi)=[x,y]_A-L^*_x\eta+L^*_y\xi;\\
{[P(x+\xi),y+\eta]_s} &=& [x-\xi,y+\eta]_s=[x,y]_A+L^*_x\eta+L^*_y\xi;\\
{[x+\xi,P(y+\eta)]_s} &=& [x+\xi,y-\eta]_s=[x,y]_A-L^*_x\eta-L^*_y\xi;\\
P[x-\xi,y-\eta]_s&=&P([x,y]_A-L^*_x\eta+L^*_y\xi)=[x,y]_A+L^*_x\eta-L^*_y\xi.
\end{eqnarray*}
Therefore, we have
$$P[x+\xi,y+\eta]_s=[P(x+\xi),y+\eta]_s+[x+\xi,P(y+\eta)]_s-P[P(x+\xi),P(y+\eta)]_s.$$
So $P$ is a paracomplex structure.\qed\vspace{3mm}

A  Lie algebroid $(A,[\cdot,\cdot]_A,a_A)$ is said to
be a {\bf quadratic Lie algebroid} if $K=\ker(a_A)$ is equipped with
a fiberwise nondegenerate $\ad$-invariant symmetric bilinear form
$(\cdot,\cdot)_+$ satisfying:
$$
a_A(x)(r, s)_+ = ([x,r]_A,s)_+ + (r,[x,s]_A)_+, \quad \forall x\in
\Gamma(A),~r,s\in\Gamma(K).
$$

\begin{defi}
A left-symmetric algebroid $(A,\cdot_A,a_A)$ is said to be a {\bf quadratic
left-symmetric algebroid} if there is a fiberwise nondegenerate
$L_x$-invariant symmetric bilinear form $(\cdot,\cdot)_+$ on $A$, i.e.
$$(x\cdot_A y,z)_++(y,x\cdot_A z)_+=a(x)(y,z)_+,\quad \forall x,y,z\in \Gamma(A).$$
Furthermore, if the symmetric bilinear form $(\cdot,\cdot)_+$ is
positive definite, we call this quadratic left-symmetric algebroid a
{\bf Riemannian left-symmetric algebroid}.
\end{defi}

Let $M$ be a Riemannian  manifold with a flat Riemannian connection
$\nabla$. Then $(TM,\nabla,\id)$ is a Riemannian left-symmetric algebroid.

\begin{pro}
Let $(A,\cdot_A,a_A)$  be a quadratic left-symmetric algebroid, such that
 \begin{equation}\label{eq:temp}
 x\cdot_A y+y\cdot_A x=0,\quad \forall~x,y\in K=\ker(a_A).
 \end{equation}
 Then $(G(A),(\cdot,\cdot)_+|_K)$ is a quadratic Lie algebroid.
\end{pro}
\pf It is not hard to see that $$a_A(x)(y,z)_+=([x,y],z)_++(y,[x,z])_+,\quad \forall x\in \Gamma(A),~y,z\in K,$$ if and only if $(y\cdot_A
x,z)_++(y,z\cdot_A x)_+=0$. On the other hand, we have $$(y\cdot_A x,z)_++(y,z\cdot_A
x)_+=-(z\cdot_A y+y\cdot_A z,x)_+.$$ Therefore, if \eqref{eq:temp} holds, then $(G(A),(\cdot,\cdot)_+|_K)$ is a quadratic Lie algebroid.\qed

\begin{defi}
Let $(A,[\cdot,\cdot]_A,a_A)$ be a real Lie algebroid, a {\bf
complex structure} on $A$ is a bundle map  $J:A\longrightarrow A$,
satisfying $J^2=-\id$ and the integrability condition:
$$J[x,y]_A=[J(x),y]_A+[x,J(y)]_A+J[J(x),J(y)]_A\quad \forall x,y\in\Gamma(A).$$
\end{defi}

\begin{pro}\label{pro:complex}
Let $(A,\cdot_A,a_A)$  be a quadratic left-symmetric algebroid, and
$\varphi:A\longrightarrow A^*$ the linear isomorphism  induced by
the nondegenerate  bilinear form $(\cdot,\cdot)_+$, i.e.
$$\langle\varphi(x),y\rangle=(x,y)_+,\quad \forall x,y\in \Gamma(A).$$
Then there exists a complex structure on the phase space $\huaP=G(A)\ltimes _{L^*}A^*$, which is  given by
$J=\left(\begin{array}{cc}
    0&-\varphi^{-1}\\\varphi&0
  \end{array}\right)$.
  More precisely, we have
\begin{eqnarray}\label{complex structure}
J(x+\xi)=-\varphi^{-1}(\xi)+\varphi(x),\quad \forall x\in
\Gamma(A),~\xi\in\Gamma(A^*).
\end{eqnarray}
\end{pro}
\pf It is obvious that $J^2=-\id$. Next, we prove the integrability
condition. Note that for all $x,y,z\in \Gamma(A)$, we have
$\varphi(x\cdot_A y)=L^*_x\varphi(y)$, which follows from
\begin{eqnarray*}
\langle\varphi(x\cdot_A y),z\rangle&=&(x\cdot_A
y,z)_+=a_A(x)(y,z)_+-(y,x\cdot_A
z)_+\\
&=&a_A(x)\langle\varphi(y),z\rangle-\langle\varphi(y),x\cdot_A
z\rangle=\langle L^*_x\varphi(y),z\rangle.
\end{eqnarray*}

 For all $x,y,z\in \Gamma(A)$, we have
 \begin{eqnarray*}
\langle J[x,y]_s,z\rangle&=&\langle \varphi(x\cdot_A y-y\cdot_A x),z\rangle=(x\cdot_A y-y\cdot_A x,z)_+;\\
\langle[J(x),y]_s,z\rangle&=&\langle[\varphi(x),y]_s,z\rangle=-\langle[y,\varphi(x)]_s,z\rangle\\
&=&-a_A(y)\langle\varphi(x),z\rangle+\langle\varphi(x),y\cdot_A z\rangle=-(y\cdot_A x,z)_+;\\
\langle[x,J(y)]_s,z\rangle&=&(x\cdot_A y,z)_+;\\
\langle J([J(x),J(y)]_s),z\rangle&=&0,
\end{eqnarray*}
which implies that
\begin{equation}\label{eq:complex1}
J[x,y]_s=[J(x),y]_s+[x,J(y)]_s+J[J(x),J(y)]_s.
\end{equation}

Similarly,   we
 have
\begin{equation}\label{eq:complex2}
J[\xi,\eta]_s=[J(\xi),\eta]_s+[\xi,J(\eta)]_s+J[J(\xi),J(\eta)]_s,\quad\forall~\xi,\eta\in \Gamma(A^*).
\end{equation}

At last, for all $x\in \Gamma(A),~\xi,\eta\in \Gamma(A^*),$ let
$y,z\in \Gamma(A)$ satisfy $\varphi(y)=\xi,~ \varphi(z)=\eta$. Then we have
\begin{eqnarray*}
\langle J[x,\xi]_s,\eta\rangle&=&\langle
-\varphi^{-1}([x,\xi]_s),\eta\rangle=\langle-\varphi^{-1}(L^*_x\varphi(y)),\eta\rangle\\
&=& -\langle x\cdot_A y,\eta\rangle=-(x\cdot_A y,z)_+;\\
 \langle[J(x),\xi]_s,\eta\rangle&=&0;\\
 \langle[x,J(\xi)]_s,\eta\rangle&=&\langle[x,-\varphi^{-1}(\xi)]_s,\eta\rangle=-(x\cdot_A y-y\cdot_A x,z)_+;\\
\langle J([J(x),J(\xi)]_s),\eta\rangle&=&\langle
J[\varphi(x),-\varphi^{-1}(\xi)]_s,\eta\rangle=\langle
J[y,\varphi(x)]_s,\eta\rangle=-(y\cdot_A x,z)_+,
\end{eqnarray*}
which implies that \begin{equation}\label{eq:complex3}
J[x,\xi]_s=[J(x),\xi]_s+[x,J(\eta)]_s+J[J(x),J(\xi)]_s.
\end{equation}

By \eqref{eq:complex1}-\eqref{eq:complex3}, $J$ is a complex structure on $\huaP$.\qed\vspace{3mm}

Let $(A,[\cdot,\cdot]_A,a_A)$ be a Lie algebroid. A {\bf complex product structure} on the Lie
algebroid $A$ is a pair $\{J,P\}$ of a complex structure $J$ and a
paracomplex structure $P$, satisfying $JP=-PJ$.

\begin{cor}
Let $(A,\cdot_A,a_A)$  be a left-symmetric algebroid. Then there exists a complex product structure $\{J,P\}$
on the phase space $\huaP=G(A)\ltimes_{L^*}A^*$, where $P$ and $J$ are given by
$(\ref{paracomplex structure})$ and $(\ref{complex structure})$
respectively.
\end{cor}

\begin{defi}
Let  $(A,[\cdot,\cdot]_A,a_A)$ be a real Lie algebroid. If there exists a complex structure
$J$ and a nondegenerate closed $2$-form $\omega$ such
that the following conditions are satisfied:
\begin{itemize}
\item[$\rm(1)$] $\omega(J(x),J(y))=\omega(x,y),$ $\quad\forall x,y\in\Gamma(A);$
\item[$\rm(2)$] $\omega(x,J(x))\neq0$,$\quad\forall x\neq 0\in\Gamma(A),$
\end{itemize}
then $\{J,\omega\}$ is called a {\bf K\"{a}hler structure} on $A$.
\end{defi}

\begin{cor}
Let $(A,\cdot_A,a_A)$  be a Riemannian left-symmetric algebroid. Then there
exists a  K\"{a}hler structure $\{J,\omega\}$ on the phase space
$\huaP=G(A)\ltimes_{L^*}A^*$, where $\omega$ is given by
$(\ref{eq:o})$ and $J$ is given by $(\ref{complex
structure})$.
\end{cor}

\section{Representations  of left-symmetric algebroids}

In this section, we study representations of left-symmetric algebroids. See \cite{cohomology of pre-Lie} for more details about representations of right-symmetric algebras.
\begin{defi}
Let $(A,\cdot_A,a_A)$ be a left-symmetric algebroid and $E$  a vector
bundle. A {\bf representation} of $A$ on $E$ consists of a pair
$(\rho,\mu)$, where $\rho:A\longrightarrow \frkD(E)$ is a representation
of $G(A)$ on $E $ and $\mu:A\longrightarrow \End(E)$ is a bundle
map, such that for all $x,y\in \Gamma(A),\ e\in\Gamma(E)$, we have
\begin{eqnarray}\label{representation condition 2}
 \rho(x)\mu(y)e-\mu(y)\rho(x)e=\mu(x\cdot_A y)e-\mu(y)\mu(x)e.
\end{eqnarray}
Denote a representation by $(E;\rho,\mu)$.
\end{defi}
Let $(A,\cdot_A,a_A)$ be a left-symmetric algebroid. It is obvious that if $(E;\rho)$ is a representation of the sub-adjacent Lie algebroid $G(A)$, then $(E;\rho,0)$ is a representation of left-symmetric algebroid $(A,\cdot_A,a_A)$.

\begin{pro}
Let $(A,\cdot_A,a_A)$ be a left-symmetric algebroid and $(E;\rho,\mu)$ its representation. Then, $(F=A\oplus E,\ast_F,a_F)$ is a left-symmetric algebroid, where
$\ast_F$ and $a_F$ are given by
\begin{eqnarray*}
(x_1+e_1)\ast_F(x_2+e_2)&=&x_1\cdot_A x_2+\rho(x_1)e_2+\mu(x_2)e_1,\\
a_F(x_1+e_1)&=&a_A(x_1),
\end{eqnarray*}
for all $x_1,x_2\in \Gamma(A),\ e_1,e_2\in \Gamma(E)$.
\end{pro}
We denote this semidirect product left-symmetric algebroid by
 $A\ltimes_{\rho,\mu} E$.

\pf Let $(E;\rho,\mu)$ be a representation. It is straightforward to  see that
$(\Gamma(F),\ast_F)$ is a left-symmetric algebra. For any $x_1,x_2\in\Gamma(A)$ and $e_1,e_2\in\Gamma(E)$, we have
\begin{eqnarray*}
(x_1+e_1)\ast_F (f(x_2+e_2))&=&x_1\cdot_A (fx_2)+\rho(x_1)e_2+\mu(x_2)e_1\\
&=& f(x_1\cdot_A x_2)+a_A(x_1)(fx_2)+f\rho(x_1) e_2+a_A(x_1)(f
e_2)+f\mu(x_2)e_1\\
&=& f((x_1+e_1)\ast_F(x_2+e_2))+a_A(x_1)(f)(x_2+e_2)\\
(f(x_1+e_1))\ast_F (x_2+e_2)&=&(fx_1)\cdot_A x_2+\rho(fx_1)e_2+\mu(x_2)(fe_1)\\
&=& f((x_1+e_1)\ast_F (x_2+e_2)).
\end{eqnarray*}
Therefore, $(F,\ast_F,a_F)$ is a left-symmetric algebroid.
\qed\vspace{3mm}

Let $(A,\cdot_A,a_A)$ be a left-symmetric algebroid, and
$(E;\rho,\mu)$ be its representation.
$\rho^*:G(A)\longrightarrow\frkD(E^*)$ is the dual representation of
$\rho$, and $\mu^*:A\longrightarrow \End(E^*)$ is defined by
$\langle \mu^*(x) \xi,e\rangle=-\langle \mu(x)e,\xi \rangle.$

\begin{pro}\label{pro:representation}
With the above notations, we have
\begin{itemize}
\item[$\rm(i)$] $(E;\rho-\mu)$ is a representation of the Lie algebroid
$G(A)$;
\item[$\rm(ii)$] $(E^*;\rho^*-\mu^*,-\mu^*)$ is a representation of $A$.
\end{itemize}
\end{pro}
\pf Since $(\Gamma(E);\rho,\mu)$ is a representation of the left-symmetric
algebra $(\Gamma(A),\cdot_A)$, $(\Gamma(E);\rho-\mu)$ is a representation of the sub-adjacent Lie algebra $(\Gamma(A),[\cdot,\cdot]_A)$. Furthermore, we have
\begin{eqnarray*}
(\rho-\mu)(f x_1)(e_1)&=&f(\rho-\mu)(x_1)(e_1),\\
(\rho-\mu)(x_1)(f e_1)&=& \rho(x_1)(f e_1)-\mu(x_1)(f e_1)\\
&=&f \rho(x_1)(e_1)+a_A(x_1)(f)e_1-f \mu(x)(e) \\
&=& f(\rho-\mu)( x_1)(e_1)+a_A(x_1)(f)e_1.
\end{eqnarray*}
Hence $(E;\rho-\mu)$ is a representation of Lie algebroid $G(A)$ on $E$. This finishes the proof of (i). (ii) can be proved similarly.
\qed

\begin{cor} With the above notations, we have
  \begin{itemize}
\item[$\rm(i)$] The left-symmetric algebroids $A\ltimes_{\rho,\mu} E$ and
$A\ltimes_{\rho-\mu,0}E$ have the same sub-adjacent
Lie algebroid $G(A)\ltimes_{\rho-\mu}E.$
\item[$\rm(ii)$] The left-symmetric algebroids
$A\ltimes_{\rho^*,0} E^*$ and $A\ltimes_{\rho^*-\mu^*,-\mu^*} E^*$
have the same sub-adjacent Lie algebroid $G(A)\ltimes_{\rho^*}E^*$.
\end{itemize}
\end{cor}

Let $(E;\rho,\mu)$ be a   representation of a left-symmetric
algebroid $(A,\cdot_A,a_A)$. In general, $(E^*;\rho^*,\mu^*) $ is not
a representation. But we have
\begin{pro}\label{dual representation condition}
Let $(A,\cdot_A,a_A)$  be a left-symmetric algebroid, and
$(E;\rho,\mu)$ be its representation. Then the following conditions
are equivalent:
\begin{itemize}
\item[$\rm(1)$]$(E;\rho-\mu,-\mu)$ is a representation of $A$.
\item[$\rm(2)$]$(E^*;\rho^*,\mu^*)$ is a representation of $A$.
\item[$\rm(3)$]$\mu(x)\mu(y)=\mu(y)\mu(x)$ for all $x,y\in \Gamma(A).$
\end{itemize}
\end{pro}

An {\bf ideal} of a left-symmetric algebroid  $(A,\cdot_A,a_A)$ is a subbundle $L\subset A$ such that $\Gamma(L)$ is an ideal of the left-symmetric algebra $(\Gamma(A),\cdot_A)$.

\begin{ex}{\rm
Let $(A,\cdot_A,a_A)$ be a transitive left-symmetric algebroid. Then
we have
\begin{itemize}
\item[$\rm(1)$] $(A;L,0)$ is a representation, and $(A^*;L^*,0)$ is a representation.
\item[$\rm(2)$] $(K;\ad,0)$ is  a representation of $A$ on $K=\ker(a_A)$, and $(K^*;{\ad}^*,0)$ is a representation of $A$ on $K^*$, where $\ad$ is the adjoint
representation of $G(A)$ on $K$.
\item[$\rm(3)$] If $K$ is an ideal of $A$, then $(K;L,R)$ is a representation of $A$ on $K$, where $R:A\longrightarrow \End(K)$ is given by $R_x(y)=y\cdot_A x,\
x\in \Gamma(A),y\in \Gamma(K)$, and $(K^*;{\ad}^*,-r^*)$ is also a
representation on $K^*$ by Proposition \ref{pro:representation}.
\end{itemize}}
\end{ex}

In general, $(A;L,R)$ is not a representation of  a left-symmetric
algebroid $(A,\cdot_A,a_A)$.
 This is different from the case of a left-symmetric algebra.\vspace{3mm}

At the end of this section, we discuss the cohomology complex for a left-symmetric algebroid $(A,\cdot_A,a_A)$ with a representation $(E;\rho,\mu)$ briefly.
  Define the set of $(n+1)$-cochains by
$$C^{n+1}(A,E)=\Gamma(\Hom(\Lambda^{n}A\otimes A,E)),\
n\geq 0.$$  For all $\omega\in C^{n}(A,E)$, define $\dM\omega:\otimes^{n+1}_{\mathbb R}\Gamma(A)\longrightarrow \Gamma(E)$ by
 \begin{eqnarray*}
 &&\dM\omega(x_1,x_2,\cdots,x_{n+1})\\
 &=&\sum_{i=1}^{n}(-1)^{i+1}\rho(x_i)\omega(x_1,x_2,\cdots,\hat{x_i},\cdots,x_{n+1})\\
 &&+\sum_{i=1}^{n}(-1)^{i+1}\mu(x_{n+1})\omega(x_1,x_2,\cdots,\hat{x_i},\cdots,x_n,x_i)\\
 &&-\sum_{i=1}^{n}(-1)^{i+1}\omega(x_1,x_2,\cdots,\hat{x_i},\cdots,x_n,x_i\cdot_A x_{n+1})\\
 &&+\sum_{1\leq i<j\leq n}(-1)^{i+j}\omega([x_i,x_j]_A,x_1,\cdots,\hat{x_i},\cdots,\hat{x_j},\cdots,x_{n+1}),
\end{eqnarray*}
for all $x_i\in \Gamma(A),i=1,2\cdots,n+1$.

\begin{pro}\label{pro:coboundary}
For all  $\omega\in C^{n}(A,E)$, we have $\dM \omega\in C^{n+1}(A,E)$. Moreover, $\dM^{2}=0$. Thus, we have a well-defined cochain complex $(C^{*+1}=\bigoplus _{n\geq 0}C^{n+1}(A,E),\dM)$.
\end{pro}

The proof is routine and we put it in the appendix.\vspace{3mm}

For $e\in\Gamma(E)$, define $\dM (e)(x)=\mu(x)e-\rho(x) e$. Then we have
$\dM^{2}(e)(x,y)=\rho(x)\rho(y)e-\rho(x\cdot_A y)e$. From this
fact, we  define the set of  $0$-cochains by
 $$C^{0}(A,E):=\{e\in
\Gamma(E)\mid \rho(x)\rho(y)e-\rho(x\cdot_A y)e=0,\ \forall
x,y\in\Gamma(A)\}.$$
Then,  we obtain a cochain complex $(C(A,E)=\bigoplus
_{n\geq 0}C^{n}(A,E),\dM)$.

\section{Deformation cohomologies of left-symmetric algebroids}

In this section, we introduce another cochain complex associated to a left-symmetric algebroid, which could
control deformations. See \cite{Marius deformation cohomology,Deformation of Lie algebroi 2} for deformations of Lie algebroids.

\begin{defi}
Let $E$ be a vector bundle over $M$, a {\bf multiderivation} of degree $n$
is a multilinear map $D\in\Hom(\Lambda^{n-1}\Gamma(E)\otimes
\Gamma(E),\Gamma(E))$, such that for all $f\in C^\infty(M)$ and sections
$x_i\in \Gamma(E)$, the following conditions are satisfied:
 \begin{eqnarray*}
D(x_1,x_2,\cdots,fx_i,\cdots,x_{n-1},x_{n})&=&f D(x_1,x_2,\cdots,x_i,\cdots,x_{n-1},x_{n}),\quad i=1,2,\cdots,n-1;\\
D(x_1,x_2,\cdots,x_{n-1},fx_{n})&=&f D(x_1,x_2,\cdots,x_{n-1},x_{n})+\sigma_D(x_1,x_2,\cdots,x_{n-1})(f)x_{n},
\end{eqnarray*}
where $\sigma_D\in \Gamma(\Hom(\Lambda^{n-1}E,TM))$ is called the
{\bf symbol}. We will denote by $\Der^n(E)$ the space of multiderivations
of $n,~n\geq 1$.
\end{defi}

 Let $(A,\cdot_A,a_A)$ be a left-symmetric algebroid. Define the set of $n$-cochains by
 $C^n(A)=\Der^n(A)$
for all $\omega\in \Der^{n}(A)$, and define the
 operator $\dM_{def}:\Der^{n}(A)\longrightarrow \Hom(\otimes^{n+1}_{\mathbb R}\Gamma(A),\Gamma(A))$,
 by\begin{eqnarray*}
 &&\dM_{def}\omega(x_1,x_2,\cdots,x_{n+1})\\
 &=&\sum_{i=1}^{n}(-1)^{i+1}x_i\cdot_A\omega(x_1,x_2,\cdots,\hat{x_i},\cdots,x_{n+1})\\
 &&+\sum_{i=1}^{n}(-1)^{i+1}\omega(x_1,x_2,\cdots,\hat{x_i},\cdots,x_n,x_i)\cdot_A x_{n+1}\\
 &&-\sum_{i=1}^{n}(-1)^{i+1}\omega(x_1,x_2,\cdots,\hat{x_i},\cdots,x_n,x_i\cdot_A x_{n+1})\\
 &&+\sum_{1\leq i<j\leq {n}}(-1)^{i+j}\omega([x_i,x_j]_A,x_1,\cdots,\hat{x_i},\cdots,\hat{x_j},\cdots,x_{n+1}),
\end{eqnarray*}
for all $\ x_i\in \Gamma(A),i=1,2\cdots,n+1$.

\begin{pro}\label{pro:deformation complex}
If $\omega\in \Der^{n}(A)$, then we have $\dM_{def}\omega\in
\Der^{n+1}(A)$, in which
$\sigma_{\dM_{def}\omega}$
is given by
 \begin{eqnarray}
\nonumber\sigma_{\dM_{def}\omega}(x_1,x_2,\cdots,x_n)
&=& \sum_{i=1}^{n}(-1)^{i+1}[a_A(x_i),\sigma_{\omega}(x_1,x_2,\cdots,\hat{x_i},\cdots,x_{n}))]_{\frkX(M)}\\
\nonumber&&+\sum_{1\leq i<j\leq n}(-1)^{i+j}\sigma_{\omega}([x_i,x_j]_A,x_1,\cdots,\hat{x_i},\cdots,\hat{x_j},\cdots,x_n)\\
&&+\sum_{i=1}^{n}(-1)^{i+1}a_A(\omega(x_1,x_2,\cdots,\hat{x_i},\cdots,x_n,x_i)).\label{eq:simbol}
\end{eqnarray}
Furthermore, we have $\dM_{def}^2=0.$ Thus, we have a cochain
complex $(C_{def}^{*+1}=\bigoplus _{n\geq 0}\Der^n(A),\dM_{def})$, which is called
the {\bf deformation complex} of $A$. The corresponding cohomology, which we denote by $H_{def}^\bullet(A)$, is called the {\bf deformation cohomology}.
\end{pro}

The proof is routine and we put it in the appendix.\vspace{3mm}

 We study deformations of
left-symmetric algebroids using the deformation cohomology.
Let $(A,\cdot_A,a_A)$ be a left-symmetric algebroid, and
$\omega\in \Der^2(A)$. Consider a $t$-parameterized family
of multiplications $\cdot_t:\Gamma(A)\otimes_{\mathbb R}\Gamma(A)\longrightarrow \Gamma(A)$ and bundle maps $a_t:A\longrightarrow TM$ given by
\begin{eqnarray*}
x\cdot_t y&=&x\cdot_A y+t\omega(x,y),\\
a_t&=&a_A+t\sigma_\omega.
\end{eqnarray*}
If  $A_t=(A,\cdot_t,a_t)$ is a left-symmetric algebroid for all $t\in I$, we say that $\omega$ generates a {\bf $1$-parameter infinitesimal
deformation} of $(A,\cdot_A,a_A)$

Since $\omega\in\Der^2(A)$, we have
$$\omega(fx,y)=f\omega(x,y),\quad \omega(x,fy)=f\omega(x,y)+\sigma_\omega(x)(f)y,$$ which implies that conditions (i) and (ii) in Definition \ref{defi:left-symmetric algebroid} are satisfied for the multiplication $\cdot_t$. Then we can deduce that  $(A_t,\cdot_t,a_t)$ is
a deformation of $(A,\cdot_A,a_A)$  if and only if
  \begin{eqnarray}\label{2-closed}
x\cdot_A\omega(y,z)-y\cdot_A\omega(x,z)+\omega(y,x)\cdot_A z-\omega(x,y)\cdot_A z    \nonumber \\
=\omega(y,x\cdot_A z)-\omega(x,y\cdot_A z)-\omega([x,y]_A,z),
\end{eqnarray}
and
\begin{eqnarray}
\omega(\omega(x,y),z)-\omega(x,\omega(y,z))&=&\omega(\omega(y,x),z)-\omega(y,\omega(x,z))\label{omega bracket}.
\end{eqnarray}
Eq. $(\ref{2-closed})$ means that $\omega$ is a $2$-cocycle, and Eq. $(\ref{omega bracket})$ means that $(A,\omega,\sigma_\omega)$ is a left-symmetric algebroid.

Two deformations $A_t$ and $A'_t$  of
a left-symmetric algebroid $A$, which are generated by $\omega$ and $\omega'$, are said to be {\bf equivalent} if
there exists a family of left-symmetric algebroid isomorphisms
$\Id+tN:A_t\longrightarrow A'_t$. A deformation is said to be
{\bf trivial} if there exists a family of left-symmetric algebroid isomorphisms
$\Id+tN:A_t\longrightarrow A$.

By straightforward computations, $A_t$ and $A'_t$  are
equivalent deformations  if and only if
\begin{eqnarray}
\omega(x,y)-\omega'(x,y)&=&x\cdot_A N(y)+N(x)\cdot_A y-N(x\cdot_A y),\label{2-exact}\\
N\omega(x,y)&=&\omega'(x,N(y))+\omega'(N(x),y)+N(x)\cdot_A N(y),\label{integral condition 1}\\
\omega'(N(x),N(y))&=&0,\label{eq:con1}\\
\sigma_\omega'\circ N&=&0,\label{eq:con2}\\
\sigma_\omega-\sigma_\omega'&=&a_A\circ
N.\label{relation of anchor def}
\end{eqnarray}

Eq. $(\ref{2-exact})$ means that $\omega-\omega'=\dM_{def}
N$.  Eqs. \eqref{eq:con1} and \eqref{eq:con2} mean that $\omega'|_{\Img(N)}=0.$ Eq. \eqref{relation of anchor def} can be obtained from Eq. \eqref{2-exact}.

We summarize the above discussion into
the following theorem:
\begin{thm}\label{thm:deformation}
Let $(A_t,\cdot_t,a_t)$ be a deformation generated by
$\omega\in\Der^2(A)$ of a left-symmetric algebroid $(A,\cdot_A,a_A)$.
Then $\omega$ is closed, i.e. $\dM_{def}\omega=0.$ Furthermore, if
two deformations $(A_t,\cdot_t,a_t)$ and $(A'_t,\cdot'_t,a'_t)$
generated by $\omega $ and $\omega'$ are equivalent,  then $\omega$
and $\omega'$ are in the same cohomology class in $H^2_{def}(A)$.
\end{thm}
Now we consider  trivial deformations of a left-symmetric algebroid
$(A,\cdot_A,a_A)$. Eqs.
$(\ref{2-exact})-(\ref{relation of anchor def})$ reduce to
\begin{eqnarray}
\omega(x,y)&=&x\cdot_A N(y)+N(x)\cdot_A y-N(x\cdot_A y),\label{Nij1}\\
N\omega(x,y)&=&N(x)\cdot_A N(y),\label{Nij2}\\
a_A\circ N&=&\sigma_\omega.\label{Nij3}
\end{eqnarray}
Again, Eq. \eqref{Nij3} can be obtained from Eq. \eqref{Nij1}. It follows from $(\ref{Nij1})$ and $(\ref{Nij2})$ that $N$ must
satisfy the following condition:
\begin{eqnarray}
N(x)\cdot_A N(y)-x\cdot_A N(y)-N(x)\cdot_A y+N^2(x\cdot_A
y)=0.\label{integral condition of Nij}
\end{eqnarray}
\begin{defi}
A bundle map
$N:A\longrightarrow A$ is called a {\bf Nijenhuis operator} on a left-symmetric algebroid $(A,\cdot_A,a_A)$  if the Nijenhuis condition
 $(\ref{integral condition of Nij})$ holds.
\end{defi}
Obviously, any Nijenhuis operator on a left-symmetric algebroid is also a
Nijenhuis operator on the corresponding  sub-adjacent Lie algebroid.

We have seen that a trivial deformation of a left-symmetric algebroid could give rise to a Nijenhuis operator. In fact, the converse is also true.

\begin{thm}\label{thm:trivial def}
  Let $N$ be a Nijenhuis operator on a left-symmetric algebroid $(A,\cdot_A,a_A)$. Then a
  deformation of  $(A,\cdot_A,a_A)$ can be
  obtained by putting
  $$
\omega(x,y)= \dM_{def}N(x,y).
  $$
  Furthermore, this deformation is trivial.
\end{thm}
\pf Since $\omega$ is exact, $\omega $ is closed naturally, i.e. Eq.
\eqref{2-closed} holds. To see that $\omega$ generates a
deformation, we only need to show that \eqref{omega bracket} holds,
which follows from the Nijenhuis condition \eqref{integral condition
of Nij}. We omit the details. At last, it is obvious that
$$
(\Id +tN)(x\cdot_t y)=(\Id+tN) (x)\cdot_A (\Id+tN)(y),\quad a_A\circ
(\Id+tN)=a_t,
$$
which implies that the deformation is trivial. \qed
\section*{Appendix}

{\bf The proof of Proposition \ref{pro:coboundary}}:

For  $\omega\in
 C^{n}(A,E)$, it is not hard to see that $\dM\omega$ is skew-symmetric with respect to the first $n$-terms.
 To see  that $\dM\omega$  is
$C^\infty(M)$-linear,
we use notations  $x\cdot
e=\rho(x)e,\ e\cdot x=\mu(x)e$
 to simplify the proof. For all $\ x_1,\cdots,x_{n+1}\in\Gamma(A),\ f\in C^\infty(M)$, we have
 \begin{eqnarray*}
&&\dM\omega(fx_1,x_2,\cdots,x_n,x_{n+1})\\
&=&(fx_1)\cdot\omega(x_2,\cdots,x_n,x_{n+1})-\sum_{i=2}^{n}(-1)^{i+1}x_i\cdot\omega(fx_1,x_2,\cdots,\hat{x_i},\cdots,x_{n+1})\\
&&+\omega(x_2,\cdots,x_n,fx_1)\cdot x_{n+1}+\sum_{i=2}^{n}(-1)^{i+1}\omega(fx_1,x_2,\cdots,\hat{x_i},\cdots,x_n,x_i)\cdot x_{n+1}\\
 &&-\omega(x_2,\cdots,x_n,(fx_1)\cdot_A{x_{n+1}})-\sum_{i=2}^{n}(-1)^{i+1}\omega(fx_1,x_2,\cdots,\hat{x_i},\cdots,x_n,x_i\cdot_A x_{n+1})\\
 &&+\sum_{j=2}^{n}(-1)^{1+j}\omega([fx_1,x_j]_A,x_2,\cdots,\hat{x_j},\cdots,x_{n+1})\\
 &&+\sum_{1\leq i<j\leq n}(-1)^{i+j}\omega([x_i,x_j]_A,fx_1,\cdots,\hat{x_i},\cdots,\hat{x_j},\cdots,x_{n+1})\\
&=&f(x_1\cdot\omega(x_2,\cdots,x_n,x_{n+1}))-f\sum_{i=2}^{n}(-1)^{i+1}x_i\cdot\omega(x_1,x_2,\cdots,\hat{x_i},\cdots,x_{n+1})\\
&&+\sum_{i=2}^{n}(-1)^{i+1}a_A(x_i)(f)\omega(x_1,x_2,\cdots,\hat{x_i},\cdots,x_{n+1})+f\omega(x_2,\cdots,x_n,x_1)\cdot x_{n+1}\\
&&+f\sum_{i=2}^{n}(-1)^{i+1}\omega(x_1,x_2,\cdots,\hat{x_i},\cdots,x_n,x_i)\cdot x_{n+1}-f\omega(x_2,\cdots,x_n,x_1\cdot_A{x_{n+1}})\\
&&-\sum_{i=2}^{n}(-1)^{i+1}f\omega(x_1,x_2,\cdots,\hat{x_i},\cdots,x_n,x_i\cdot_A x_{n+1})\\
&&+f\sum_{j=2}^{n}(-1)^{1+j}\omega([x_1,x_j]_A,x_2,\cdots,\hat{x_j},\cdots,x_{n+1})\\
&&-\sum_{j=2}^{n}(-1)^{1+j}a_A(x_j)(f)\omega(x_1,x_2,\cdots,\hat{x_j},\cdots,x_{n+1})\\
&&+f\sum_{1\leq i<j\leq n}(-1)^{i+j}\omega([x_i,x_j]_A,x_1,\cdots,\hat{x_i},\cdots,\hat{x_j},\cdots,x_{n+1})\\
&=&f\dM\omega(x_1,x_2,\cdots,x_n,x_{n+1}).
\end{eqnarray*}
  Since $\dM\omega$ is skew-symmetric with respect to the first $n$-terms,
we deduce that $\dM\omega$ is $C^\infty(M)$-linear with respect to the first $n$
terms. Similarly, we can prove that $\dM\omega$ is $C^\infty(M)$-linear with respect to the last term.

Furthermore, since the coboundary operator $\dM$ is also a
left-symmetric algebra coboundary operator, we have $\dM^{2}=0$.
This finishes the proof. \qed\vspace{3mm}

{\bf The proof of Proposition \ref{pro:deformation complex}}:

Obviously, for $\omega\in \Der^{n}(A)$, we have $\dM_{def}\omega\in\Hom(\Lambda^n\Gamma(A)\otimes_{\mathbb R} \Gamma(A),\Gamma(A))$. For $n\geq3$, we have
 \begin{eqnarray*}
&&\dM_{def}\omega(fx_1,x_2,\cdots,x_n,x_{n+1})\\
&=&(fx_1)\cdot_A\omega(x_2,\cdots,x_n,x_{n+1})+\sum_{i=2}^{n}(-1)^{i+1}x_i\cdot_A\omega(fx_1,x_2,\cdots,\hat{x_i},\cdots,x_{n+1})\\
&&+\omega(x_2,\cdots,x_n,fx_1)\cdot_A x_{n+1}+\sum_{i=2}^{n}(-1)^{i+1}\omega(fx_1,x_2,\cdots,\hat{x_i},\cdots,x_n,x_i)\cdot_A x_{n+1}\\
 &&-\omega(x_2,\cdots,x_n,(fx_1)\cdot_A{x_{n+1}})-\sum_{i=2}^{n}(-1)^{i+1}\omega(fx_1,x_2,\cdots,\hat{x_i},\cdots,x_n,x_i\cdot_A x_{n+1})\\
 &&+\sum_{j=2}^{n}(-1)^{1+j}\omega([fx_1,x_j]_A,x_2,\cdots,\hat{x_j},\cdots,x_{n+1})\\
 &&+\sum_{1\leq i<j\leq n}(-1)^{i+j}\omega([x_i,x_j]_A,fx_1,\cdots,\hat{x_i},\cdots,\hat{x_j},\cdots,x_{n+1})\\
&=&f(x_1\cdot_A\omega(x_2,\cdots,x_n,x_{n+1}))-f\sum_{i=2}^{n}(-1)^{i+1}x_i\cdot_A\omega(x_1,x_2,\cdots,\hat{x_i},\cdots,x_{n+1})\\
&&+\sum_{i=2}^{n}(-1)^{i+1}a_A(x_i)(f)\omega(x_1,x_2,\cdots,\hat{x_i},\cdots,x_{n+1})\\
&&+f\omega(x_2,\cdots,x_n,x_1)\cdot_A x_{n+1}+\sigma_{\omega}(x_2,\cdots,x_n)(f)x_1\cdot_A x_{n+1}\\
&&+f\sum_{i=2}^{n}(-1)^{i+1}\omega(x_1,x_2,\cdots,\hat{x_i},\cdots,x_n,x_i)\cdot_A x_{n+1}\\
&&-f\omega(x_2,\cdots,x_n,x_1\cdot_A{x_{n+1}})-\sigma_{\omega}(x_2,\cdots,x_n)(f)x_1\cdot_A x_{n+1}\\
&&-\sum_{i=2}^{n}(-1)^{i+1}f\omega(x_1,x_2,\cdots,\hat{x_i},\cdots,x_n,x_i\cdot_A x_{n+1})\\
&&+f\sum_{j=2}^{n}(-1)^{1+j}\omega([fx_1,x_j]_A,x_2,\cdots,\hat{x_j},\cdots,x_{n+1})\\
&&-\sum_{j=2}^{n}(-1)^{1+j}a_A(x_j)(f)\omega(x_1,x_2,\cdots,\hat{x_j},\cdots,x_{n+1})\\
&&+f\sum_{1\leq i<j\leq n}(-1)^{i+j}\omega([x_i,x_j]_A,x_1,\cdots,\hat{x_i},\cdots,\hat{x_j},\cdots,x_{n+1})\\
&=&f\dM_{def}\omega(x_1,x_2,\cdots,x_n,x_{n+1}),
\end{eqnarray*}
which implies that $\dM_{def}\omega$ is $\CWM$-linear in the first component. Since $\dM_{def}\omega$ is skew-symmetric in the first $n$ components, we deduce that  $\dM_{def}\omega$ is $\CWM$-linear in the first $n$ components. One can verify that this fact also holds for the case of $n=1,2$, we omit details.

As for the last component,  we  have
 \begin{eqnarray*}
 &&\dM_{def}\omega(x_1,x_2,\cdots,x_n,fx_{n+1})\\
 &=&\sum_{i=1}^{n}(-1)^{i+1}x_i\cdot_A\omega(x_1,x_2,\cdots,\hat{x_i},\cdots,fx_{n+1})\\
 &&+\sum_{i=1}^{n}(-1)^{i+1}\omega(x_1,x_2,\cdots,\hat{x_i},\cdots,x_n,x_i)\cdot_A (fx_{n+1})\\
 &&-\sum_{i=1}^{n}(-1)^{i+1}\omega(x_1,x_2,\cdots,\hat{x_i},\cdots,x_n,x_i\cdot_A (fx_{n+1}))\\
 &&+\sum_{1\leq i<j\leq n}(-1)^{i+j}\omega([x_i,x_j]_A,x_1,\cdots,\hat{x_i},\cdots,\hat{x_j},\cdots,fx_{n+1})\\
&=&\sum_{i=1}^{n}(-1)^{i+1}x_i\cdot_A(f\omega(x_1,x_2,\cdots,\hat{x_i},\cdots,x_{n+1})+\sigma_{\omega}(x_1,x_2,\cdots,\hat{x_i},\cdots,x_{n})(f)x_{n+1})\\
&&+f\sum_{i=1}^{n}(-1)^{i+1}\omega(x_1,x_2,\cdots,\hat{x_i},\cdots,x_n,x_i)\cdot_A x_{n+1}\\
&&+\sum_{i=1}^{n}(-1)^{i+1}a_A(\omega(x_1,x_2,\cdots,\hat{x_i},\cdots,x_n,x_i))(f)x_{n+1}\\
&&-\sum_{i=1}^{n}(-1)^{i+1}\omega(x_1,x_2,\cdots,\hat{x_i},\cdots,x_n,f(x_i\cdot_A x_{n+1})+a_A(x_i)(f)x_{n+1})\\
&&+f\sum_{1\leq i<j\leq n}(-1)^{i+j}\omega([x_i,x_j]_A,x_1,\cdots,\hat{x_i},\cdots,\hat{x_j},\cdots,x_{n+1})\\
&&+\sum_{1\leq i<j\leq n}(-1)^{i+j}\sigma_{\omega}([x_i,x_j]_A,x_1,\cdots,\hat{x_i},\cdots,\hat{x_j},\cdots,x_n)(f)x_{n+1}\\
&=&f\dM_{def}\omega(x_1,x_2,\cdots,x_n,x_{n+1})+\sum_{i=1}^{n}(-1)^{i+1}[a_A(x_i),\sigma_{\omega}(x_1,x_2,\cdots,\hat{x_i},\cdots,x_{n}))]_{\frkX(M)}(f)x_{n+1},\\
&&+\sum_{1\leq i<j\leq n}(-1)^{i+j}\sigma_{\omega}([x_i,x_j]_A,x_1,\cdots,\hat{x_i},\cdots,\hat{x_j},\cdots,x_n)(f)x_{n+1}\\
&&+\sum_{i=1}^{n}(-1)^{i+1}a_A(\omega(x_1,x_2,\cdots,\hat{x_i},\cdots,x_n,x_i))(f) x_{n+1},
\end{eqnarray*}
which implies that $\dM_{def}\omega\in \Der^{n+1}(A)$ and  $\sigma_{\dM\omega}$ is given by \eqref{eq:simbol}.

At last, the conclusion $\dM_{def}^2=0$ can be obtained similarly as the case of left-symmetric algebras. We omit the details. The proof is completed. \qed

\end{document}